\documentstyle[twoside,12pt]{article}
\def \BB{{\cal B}}
\def\bee{\begin{equation}}
\def\ee{\end{equation}}
\input epsf

\begin{document}

\thispagestyle{empty}
\bigskip\bigskip
\centerline{    }
\vspace*{2.5cm}

\begin{center}
{\Large \bf Numerical evidence in favor of the Arenstorf formula}\\*[7mm]

{\large \sl Marek Wolf}\\*[5mm]

Institute of Theoretical Physics, University of Wroc{\l}aw\\
Pl.Maxa Borna 9, PL-50-204 Wroc{\l}aw, Poland, e-mail:
mwolf@ift.uni.wroc.pl\\
\end{center}

\begin{abstract}

The formula $\lim_{N\rightarrow\infty}\sum_{p<N,p,p+2~both~prime}
\log(p)\log(p+2) = C_2$ is tested on the computer up to
$N=2^{40}\approx 1.1\times 10^{12}$  and very good agreement  is
found.

\end{abstract}

\vspace*{2cm}

The recent paper ``There Are Infinitely Many Prime Twins'' by R.F.
Arenstorf \cite{Arenstorf} has raised a lot of excitement. The
author claims to proved that:

\bee
\lim_{N\rightarrow\infty}\sum_{p<N,p,p+2~ twins}\log(p)\log(p+2) = C_2
\label{main}
\ee
where the twin  constant
\bee
C_2 \equiv
2\prod_{p > 2} \biggl( 1 - {1 \over (p - 1)^2}\biggr)
=1,32032363169373914785562\ldots.
\ee

Waiting for the formal approval of this result by mathematical
community I have run the computer program to check the validity of (\ref{main}). Even Hardy
and Littlewood in their famous paper \cite{Hardy} have presented tables with numerical
verification of their conjectures based on the existing that time data up to 9 000 000.
Because  the program I have written several years ago
performs all operations on bits 
it was natural to store the data representing the
actual arithmetical mean value of $\log(p)\log(p+2)$  at values of $N$
forming the geometrical
progression with the ratio 2, i.e. at $N=2^{22}, 2^{23}, \ldots, 2^{39},
2^{40}$. The results are presented in the Table 1.

\newpage

\vskip 0.4cm
\begin{center}
{\sf TABLE {\bf I}}\\

\bigskip
\bigskip

\begin{tabular}{|c||c||c||c|} \hline
$ N $ & $1/N\sum_{p<N}\log(p)\log(p+2) $ & $ 1/N\sum_{p<N}\log(p)\log(p+2)/C_2 $ \\ \hline
$2^{22}=$        4194304 &   1.330875543&      1.00799191245 \\ \hline
$2^{23}=$        8388608 &   1.325154123&      1.00365856579 \\ \hline
$2^{24}=$       16777216 &   1.323501313&      1.00240674443 \\ \hline
$2^{25}=$       33554432 &   1.320938577&      1.00046575325 \\ \hline
$2^{26}=$       67108864 &   1.319310330&      0.99923253570 \\ \hline
$2^{27}=$      134217728 &   1.320265943&      0.99995630684 \\ \hline
$2^{28}=$      268435456 &   1.319095515&      0.99906983694 \\ \hline
$2^{29}=$      536870912 &   1.319679380&      0.99951205023 \\ \hline
$2^{30}=$     1073741824 &   1.320096901&      0.99982827612 \\ \hline
$2^{31}=$     2147483648 &   1.320047000&      0.99979048210 \\ \hline
$2^{32}=$     4294967296 &   1.320350510&      1.00002035747 \\ \hline
$2^{33}=$     8589934592 &   1.320423613&      1.00007572498 \\ \hline
$2^{34}=$    17179869184 &   1.320490713&      1.00012654539 \\ \hline
$2^{35}=$    34359738368 &   1.320290443&      0.99997486310 \\ \hline
$2^{36}=$    68719476736 &   1.320309503&      0.99998929874 \\ \hline
$2^{37}=$   137438953472 &   1.320365187&      1.00003147324 \\ \hline
$2^{38}=$   274877906944 &   1.320351882&      1.00002139677 \\ \hline
$2^{39}=$   549755813888 &   1.320340769&      1.00001297956 \\ \hline
$2^{40}=$  1099511627776 &   1.320322532&      0.99999916675 \\ \hline
\end{tabular} \\
\end{center}
\vskip 0.4cm

As it can be seen from above table there is no apparent dependence on $N$.
Indeed, trying to find heuristically the dependence on $N$ we can argue that the
probability to find twin pair around $x$ is $C_2/\log^2(x)$ and hence the mean
expectation value of the product $\log(p)\log(p+2)$ for $p$ and $p+2$ on both sides of
$x$ ($p=x-1, p+2=x+1$) does not depend on $x$  and we have simply that

\bee
\sum_{p<N,p,p+2 ~twins}\log(p)\log(p+2) = C_2N
\ee

It can be contrasted with the calculation of the Brun constant

\bee
\BB_2=\left( {1 \over 3} + {1\over 5}\right) + \left( {1 \over 5} + {1\over
7}\right) + \left( {1 \over 11} + {1\over 13}\right) + \ldots < \infty.
\label{def_B2}
\ee

The probability
to find a pair of twins in the vicinity of $x$ is $2C_2/\log^2(x)$, so
the expected value of the finite approximation to the Brun constant can be
estimated as follows:
\bee
\BB_2(x) = \BB_2(\infty) - \sum_{p~{\rm twin~prime}>x} {1 \over p} \approx
\BB_2 - 4 c_2\int_x^\infty {du \over u\log^2(u)} = \BB_2 - {4c_2\over
\log(x)}.
\label{B2}
\ee
It means that the plot of finite approximations $\BB_2(x)$ to the
original Brun
constant is a linear function of $1/\log(x)$ \cite{Brun} and from the partial sum $\BB(x)$
calculated on the computer up to $x$ the limiting value can be extrapolated by adding
$4C_2/\log(x)$: $\BB=\BB(x)+4C_2/\log(x)$.
To gain some idea what value of the limit can be the extrapolated from
numbers in Table I the Figure 1 presents actual values of the mean value of $\log(p)\log(p+2)$
plotted against $1/N$. Fitting the straight line to these points by least square method
gives the intercept (what corresponds to $N=\infty$) 1.3200385787619. In fact we see in Table I
shortage of twins in the interval $(2^{26}, 2^{31})$ and in the next intervals
some surplus of twins. Thus
skipping the first 10 points and fitting straight line in the interval $(2^{32}, 2^{40})$
(in fact only two points are needed to determine straight line!) I got for the limiting value
of the intercept 1.3203501777.

\newpage

\thispagestyle{empty}

\begin{figure}[p]
\centering
\vspace{-1.7cm}
\hspace{7cm}
\epsfxsize=14.5truecm
\epsfbox{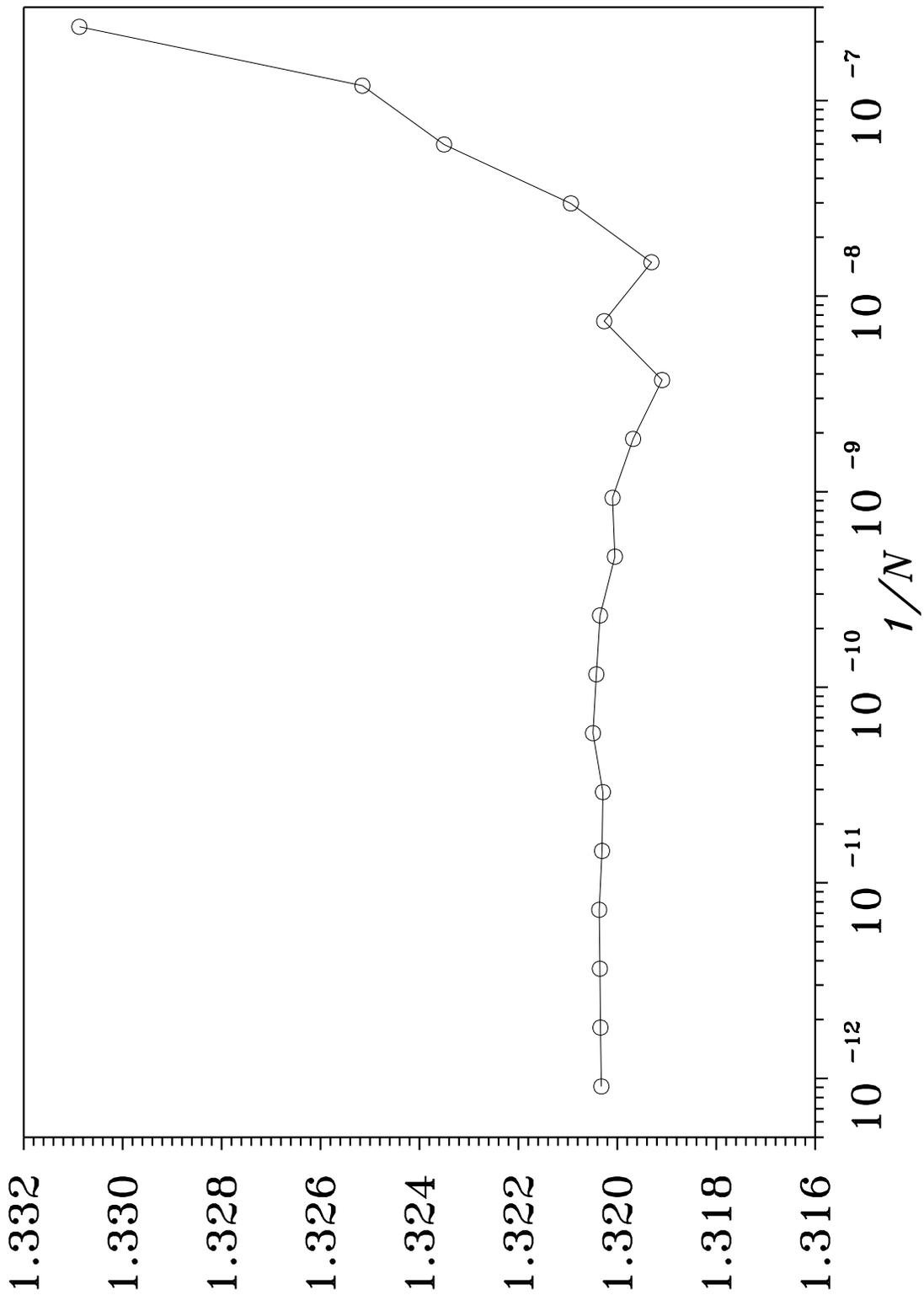}
\vspace{0.7cm}
\caption{The plot of the actual mean values of $\log(p)\log(p+2)$ against $1/N$. Notice that
on  the $x$ axis there is a logarithmic scale.}
\end{figure}

\end{document}